\documentclass[12pt]{amsart}
\usepackage{amscd,amsmath,amssymb,amsfonts}
\usepackage[cmtip, all]{xy}
\theoremstyle{plain}
\newtheorem{thm}{Theorem}

\newtheorem{prop}[thm]{Proposition}

\theoremstyle{definition}
\newtheorem{defn}[thm]{Definition}

\newtheorem{rmk}[thm]{Remark}

\numberwithin{thm}{section}
\numberwithin{equation}{section}

\newcommand{\ga}[2]{\begin{gather}\label{#1}#2 \end{gather}}
\newcommand{\gr}{{\rm gr}}

\newcommand{\surj}{\twoheadrightarrow}

\newcommand{\sI}{{\mathcal I}}

\newcommand{\sO}{{\mathcal O}}

\newcommand{\sW}{{\mathcal W}}

\newcommand{\C}{{\mathbb C}}

\newcommand{\F}{{\mathbb F}}

\renewcommand{\P}{{\mathbb P}}
\newcommand{\Q}{{\mathbb Q}}

\newcommand{\Z}{{\mathbb Z}}

\begin{document}

\title[Rational Singularities and Rational Points]{
Rational Singularities and Rational Points
}
\author{Manuel Blickle}
\email{manuel.blickle@uni-essen.de}
\author{H\'el\`ene Esnault}
\email{esnault@uni-essen.de}
\address{
Universit\"at Duisburg-Essen, Mathematik, 45117 Essen, Germany}
\thanks{Partially supported by the
DFG-Schwerpunkt ``Komplexe Mannigfaltigkeiten'' and by the DFG Leibniz Preis.}

\date{January 5, 2006}
\begin{abstract}
If $X$ is a projective, geometrically irreducible variety defined over a finite field $\F_q$, such that it is smooth and its Chow group of 0-cycles fulfills base change, i.e. $CH_0(X\times_{\F_q}\overline{ \F_q(X)})=\Q$, then the main theorem of \cite{EFano} (Corollary 1.2) asserts that its number of rational points satisfies
$|X(\F_q)| \equiv 1$ modulo $q$. If $X$ is not smooth, this is no longer true. Indeed J. Koll\'ar constructed an example of a rationally connected surface over $\F_q$ without any rational points. We show how Theorem 1.1 of \cite{BBE} indicates how to define a notion of Witt-rational singularities in characteristic $p>0$ (Definition \ref{defn2.3}). Then if $X/\F_q$ is  a projective, geometrically irreducible variety, such that it has Witt-rational singularities  and its Chow group of 0-cycles fulfills base change, then
$|X(\F_q)| \equiv 1$ modulo $q$ (Theorem \ref{thm1.1}).

\end{abstract}
\maketitle
\begin{quote}

\end{quote}

\section{Introduction}
If $X$ is a projective, geometrically irreducible variety defined over a finite field $\F_q$, such that it is smooth and its Chow group of 0-cycles fulfills base change, i.e. $CH_0(X\times_{\F_q}\overline{ \F_q(X)})=\Q$, then the main theorem of \cite{EFano} (Corollary 1.2) asserts that its number of rational points fulfills
$|X(\F_q)| \equiv 1$ modulo $q$. Smoothness enters into the proof at two levels. One uses that the Chow group $CH_*(X\times X)$ acts on cohomology $H^*(X)$ via correspondences, where $H^i(X)$ denotes either crystalline cohomology or $\ell$-adic cohomology. Then,  one uses purity in order to show that
on the first coniveau $N^1H^*(X)$ of cohomology, the eigenvalues of the geometric Frobenius are divisible by $q$ as algebraic integers if $H^*(X)$ is $\ell$-adic cohomology, or have slopes $<1$ if $H^*(X)$ is crystalline cohomology.
Once this has been established, the argument is simple. Bloch's decomposition of the diagonal, which is equivalent to the Chow group condition, forces cohomology to have coniveau 1, which in turn implies the result using Grothendieck-Lefschetz trace formula. If $X$ is no longer smooth, $CH_*(X\times X)$ does not act on cohomology, one does not have purity, and it may happen that $H^i(X)=N^1H^i(X)$
and yet the eigenvalues are not divisible, as demonstrated by the example of a rational nodal curve, for which $N^1H^1(X)=H^1(X)=\Q_\ell(0).$ Indeed, J. Koll\'ar constructed an example
of a rationally connected surface over $\F_q$ without any rational points (see section \ref{sec3.3}).

On the other hand, we show in \cite{BBE}, Theorem 1.1 that if $X$ is proper over a perfect field $k$ of characteritic $p>0$, then the slope $<1$ piece of rigid cohomology $H^i(X/K)^{<1}$ is computed by Witt vector cohomology
$H^i(X, W\sO_X)_K$, where $K={\rm Frac}(W(k)$ is the field of fractions of the ring of Witt vectors $W(k)$ and $_K$ means $\otimes_{W(k)}K$.  We explain in section \ref{sec2} how this suggests a good definition of rational singularities over a perfect field of characteristic $p>0$, which we call Witt-rational singularities.
Indeed, if $X$ is proper over $k$ of characteristic 0, the corner piece
${\rm gr}^0_FH^i(X)$  of its $F$-filtration on de  Rham cohomology receives surjectively
$\sO_X$-cohomology $H^i(X, \sO_X)$ (\cite{E}, Proposition 1.2). Thus
if $X$ has rational singularities, that is if
\ga{1}{\sigma^*: \sO_X\xrightarrow{\cong} R\sigma_*\sO_Y} for a desingularization $\sigma: Y\to X$, then $\sigma^*$ induces an isomorphism
\ga{2}{\sigma^*: {\rm gr}^0_FH^i(X) \to {\rm gr}^0_F H^i(Y).}
The analogy between the corner piece ${\rm gr}^0_FH^i(X)$ in de Rham cohomology in characteristic 0 with
the slope $<1$ piece of rigid cohomology over a perfect field of characteristic $p>0$ for a proper variety suggests that a good definition of rational singularities
in characteristic $p>0$ should imply \eqref{2} with ${\gr}^0_F H^i(X)$ replaced with $H^i(X/K)^{<0}$. Since in general we do not have a desingularization at disposal, it should be replaced by an alteration $\sigma: Y\to X$ which is generically \'etale. Then
\ga{3}{\sigma^*: H^i(X/K)^{<1}\to H^i(Y/K)^{<1} \  {\rm should \ be \ injective}.}
So we define  $X$ over a perfect field $k$ of characteristic $p>0$
to have {\it Witt-rational singularities} (see Definition \ref{defn2.3}) if for a generically \'etale alteration $\sigma: Y\to X$,
\ga{4}{\sigma^*: (W\sO_X)_K\to R\sigma_*(W\sO_Y)_K \ \  {\rm splits}.}
This definition does not depend on the choice of $\sigma$.
In characteristic 0, as observed by Kov\'acs (\cite{Ko}, Theorem 1), using Serre duality and Grauert-Riemenschneider vanishing theorem, one can replace the characterization \eqref{1} by $\sigma^*: \sO_X\to R\sigma_*\sO_Y$ being split for $\sigma: Y\to X$ an alteration.
Furthermore, our definition forces
\eqref{3} to be fulfilled.

Witt-rational singularities are  a broader class than the more restrictive class of rational singularities defined by Kempf \cite{Ke}. Indeed, his class is defined only under the condition that the singularity admits a resolution, and then the request consists of both \eqref{1} and Grauert-Riemenschneider vanishing theorem.
For those singularities for which one knows the existence of a resolution, his definition is more restrictive (see the example discussed in section \ref{sec2.3}).

Our main theorem, proven in section 3, is then
\begin{thm} \label{thm1.1}
Let $X$ be a geometrically irreducible projective variety defined over the finite field $\F_q$. If $X$ has Witt-rational singularities and its Chow group of 0-cycles fulfills base change $CH_0(X\times_{\F_q}
\overline{\F_q(X)})=\Q$, then $|X(\F_q)| \equiv 1$ modulo $q$.
\end{thm}
Theorem \ref{thm1.1} partly answers the questions raised in \cite{EFano}, section 3 and \cite{Euneq}, 5.2.

The method follows the pattern of the proof of Corollary 1.2 in
\cite{EFano}. Bloch's decomposition of the diagonal holds true even for singular varieties. The crucial observation is that even if $CH_*(X\times X)$ does not act on cohomology, the decomposition of the diagonal implies that the piece of cohomology of $X$ which survives in the cohomology of an alteration lies indeed in the first coniveau level:
\begin{prop}\label{prop1.2}
Let $X$ be  a projective variety defined over a field $k$ of
finite type over $\F_p$ such that $CH_0(X\times_k \overline{k(X)})=\Q$. Then if  $\sigma: Y\to X$ is a generically \'etale alteration, one has
\ga{3.1}{\sigma^*H^i(X)\subset N^1H^i(Y).
}
Here $H^i(X)$ denotes either $\ell$-adic cohomology $H^i(\bar{X}, \Q_\ell)$ or
rigid cohomology $H^i(\bar{X}/{\rm Frac}(W(\bar{k})))$.
\end{prop}
{\it Acknowledgments:} The notion of Witt-rational singularities defined in this note relies on \cite{BBE}. We thank P. Berthelot and S. Bloch for the interesting discussions we had on Witt vector cohomology while thinking of \cite{BBE}. We thank N. Fakhruddin and E. Viehweg for discussions on rational singularities in characteristic 0. Finally we thank J. Koll\'ar for sending us his inspiring example and for discussions.
\section{Rational singularites} \label{sec2}
\subsection{Characteristic 0}
Let $X$ be a variety defined over a field $k$ of characteristic 0. Then to say that $X$ has rational singularities is to say that if $\sigma: Y\to X$ is a desingularization, then
\ga{2.1}{\sigma^*: \sO_X\to R\sigma_*\sO_Y} is an isomorphism. Grauert-Riemenschneider vanishing theorem asserting that $R^i\sigma_*\omega_Y=0$ for $i>0$ implies that condition \eqref{2.1} is equivalent to saying that
\ga{2.2}{X \ {\rm Cohen-Macaulay \ and} \ \sigma^*: \omega_X\to \sigma_*\omega_Y \ {\rm is \ an \ isomorphism}.}
Kov\'acs \cite{Ko}, Theorem 1, observes then that the birationaliy of $\sigma$ in formulation \eqref{2.2} is not important, for example if $\sigma: Y\to X$ was instead an alteration, that is proper generically finite with $Y$ smooth, then \eqref{2.2} would be equivalent to the conditions $X$ being Cohen-Macaulay and $\sigma^*: \omega_X\to \sigma_*\omega_Y $  being split. So \cite{Ko}, Theorem 1, applied to alterations says
\begin{prop} \label{prop2.1}
Let $X$ be a variety over a field of characteristic 0. Then $X$ has rational singularities if and only if for some (and then any) alteration $\sigma: Y\to X$, the natural map $\sigma^*: \sO_X\to R\sigma_*\sO_Y$ splits in the bounded derived category of coherent sheaves on $X$.
\end{prop}
One of the important properties of projective varieties with rational singularities is that its $\sO$-cohomology computes the corner piece ${\rm gr}^0_F$ of the $F$-Hodge filtration on de Rham cohomology:
\begin{prop} \label{prop2.2}
Let $X$ be a projective variety with rational singularities. Then
there is a functorial isomorphism $H^i(X, \sO_X)\to {\rm gr}^0_FH^i(X)$ for all $i\ge 0$. If $\sigma: Y\to X$  is a desingularization, then $\sigma^*:
{\rm gr}^0_F H^i(X)\to {\rm gr}^0_F H^i(Y)$ is an isomorphism and if $\sigma: Y\to X$ is an alteration, then $\sigma^*:
{\rm gr}^0_F H^i(X)\to {\rm gr}^0_F H^i(Y)$ is injective.
\end{prop}
\begin{proof}
By \cite{E}, Proposition 1.2, there is a functorial surjection $H^i(X, \sO_X)\to {\rm gr}^0_FH^i(X)$. On the other hand, let $\sigma: Y\to X$ be a desingularization. Then one has the commutative diagram
\ga{2.3}{\begin{CD}
H^i(Y, \sO_Y)@>\cong >> {\rm gr}^0_FH^i(Y)\\
@A\pi^* AA @A\pi^* AA\\
 H^i(X, \sO_X)@>{\rm surj} >> {\rm gr}^0_FH^i(X)
         \end{CD}}
By the rationality condition, $\pi^*$ on the left is an isomorphism. This implies that $\pi^*$ on the right and the surjection are  isomorphisms as well.
If now $\sigma: Y\to X$ is an alteration, then Proposition \ref{prop2.1}
together with the case of a desingularization in Proposition \ref{prop2.2} allow to conclude that $\sigma^*: {\rm gr}^0_FH^i(X)\to {\rm gr}^0_FH^i(Y)$ is injective.
\end{proof}
\subsection{Characteristic $p>0$}
Let $X$ be a variety defined over a perfect field $k$ of characteristic $p>0$.
Then one has at disposal Berthelot's rigid cohmology $H^i(X/K)$, together with its slope filtration. Here $K={\rm Frac}(W(k))$ is the field of fractions of the ring of Witt vectors $W(k)$ of $k$. The pair $(H^i(X/K), {\rm slope \ filtration})$ behaves analogously to the pair (de Rham cohomology, F-filtration) in characteristic 0. We denote by $H^i(X/K)^{<1}$ the slope $<1$ piece.
Thus, a definition of rational singularities in characteristic $p>0$ should imply injectivity of $\sigma^*: H^i(X/K)^{<1}\to H^i(Y/K)^{<1}$ where $\sigma: Y\to X$ is a generically \'etale alteration, if $X$ is projective.
By \cite{deJ},  Theorem 4.1 and Remark 4.2, they exists, after a finite separable extension $k'\supset k$, a generically \'etale alteration $\sigma: Y\to X$.
Since the slope decomposition is defined over $k$ (see \cite{BBE}, Section 5.2), we may replace
$k$ by $k'$ and assume that $\sigma$ is defined over $k$.
Moreover, by \cite{BBE}, Theorem 1.1 one has
\ga{2.4}{X  \ {\rm proper} \Longrightarrow H^i(X/K)^{<1}\cong H^i(X, W\sO_X)_K\cong H^i(X, W\sO_{X,K})}
where the identification of the slope $<1$ piece of rigid cohomology with Witt vector cohomology is functorial. Here $W\sO_{X,K}=(W\sO_X)\otimes_{W(k)}K$. This generalizes to the non-smooth case the classical slope theorem of Bloch-Illusie (see Introduction of \cite{BBE} and references there).

As in \cite{BBE}, section 2, we consider for an alteration $\sigma: Y\to X$ which is generically \'etale the sheaves of abelian groups $R^i\sigma_* W\sO_{Y,K}$ on $X$.
In the bounded derived category
of sheaves of abelian groups on $X$ taken with the Zariski topology, one has the map
\ga{2.5}{\sigma^*: W\sO_{X,K}\to R\sigma_* W\sO_{Y, K}.}
\begin{defn} \label{defn2.3}
Let $X$ be a variety defined over a perfect field $k$ of characteristic $p>0$. Then $X$ has {\it Witt-rational singularities} if there is a finite separable extension of $k'\supset k$ and a generically \'etale alteration $\sigma: Y\to X$ over $k'$ for which \eqref{2.5} splits in the derived category of sheaves of
abelian groups on $X$ taken with the Zariski topology.
\end{defn}
As it stands, the definition might a priori depend on the choice of $\sigma$. To show independency, one would have to show that \eqref{2.5} splits assuming $X$ smooth and $\sigma: Y\to X$ any generically \'etale alteration. In characteristic 0, to show the corresponding statement with $W\sO_{X,K}$ replaced by $\sO_X$, one would use duality and Grauert-Riemenschneider vanishing theorem, which is analytic. We do not address  this delicate point here and instead remark the following.
\begin{prop} \label{prop2.4}
Let $X$ be a proper variety defined over a perfect field $k$ of characteristic $p>0$. Then if $X$ has Witt-rational singularities and $\sigma: Y\to X$ is any alteration, then $\sigma^*: H^i(X/K)^{<1}\to H^i(Y/K)^{<1}$ is injective.
\end{prop}
\begin{proof}
Let $\sigma_0: Y_0\to X$ be an alteration with \eqref{2.5} splitting
for $\sigma_0$ and let $\sigma: Y\to X$ be any alteration. One constructs
\ga{2.6}{\begin{CD}
Z@>\tau>> Y_0\\
@V\tau_0 VV @V\sigma_0 VV \\
Y@>\sigma >> X
         \end{CD}}
with $\tau, \tau_0$ alterations. By \cite{deJ},  Theorem 4.1 and Remark 4.2, they exist after a finite separable extension $k'\supset k$. Since the slope decomposition is defined over $k$ (see \cite{BBE}, Section 5.2), we may replace
$k$ by $k'$ and assume $\tau, \tau_0$ are defined over $k$. On the other hand,
by \cite{BBE}, Theorem 1.1, one has $H^i(X/K)^{<1}=H^i(X, W\sO_{X,K})$. By the Witt-rationality assumption, one has $\sigma_0^*: H^i(X, W\sO_{X,K})
\to H^i(Y_0, W\sO_{Y_0, K})$ injective. By the classical Bloch-Illusie slope theorem, one has
$H^i(Y_0, W\sO_{Y_0, K})=H^i(Y_0/K)^{<1}$.  By smoothness of $Z$ and $Y_0$,  $\tau^*: H^i(Y_0/K)\to H^i(Z/K)$ is injective on the whole crystalline cohomology (as one has duality and a trace map), a fortiori on the slope $<1$ piece. Thus
\ga{2.7}{\tau^*\circ \sigma_0^*=\tau_0^*\circ \sigma^* : H^i(X/K)^{<1}\to H^i(Z/K)^{<1} \ {\rm injective}.}
As
$\sigma^*: H^i(X/K)\to H^i(Y/K)$ factors through $\tau_0^*\circ \sigma^*$, it is injective as well. This shows the proposition.
\end{proof}
\begin{rmk}\label{rmk2.5}
Let $X$ be a variety defined over a perfect field of characteristic $p>0$. Assume $X$ admits a resolution of singularities $\sigma: Y\to X$ with the property that the natural map $\sigma^*: \sO_X\to R\sigma_*\sO_Y$
be an isomorphism. Then $X$ has Witt-rational singularities.
\end{rmk}
\begin{proof}
Indeed, one has $R^i\sigma_*\sO_Y=0$ for all $i\ge 1$. Consequently,
$R^i\sigma_*W\sO_Y=0$ for all $i\ge 1$, thus a fortiori $R^i\sigma_*W\sO_{Y, K}=0$
for all $i\ge 1$. On the other hand, the extra condition $\sigma_*\sO_Y=\sO_X$
implies $\sigma_*W\sO_Y=W\sO_X$. Thus the splitting in Definition \ref{defn2.3} is
in fact an isomorphism.
\end{proof}

\subsection{Example} \label{sec2.3}
It is easy to produce examples for which the converse of Remark \ref{rmk2.5} is not true. Indeed, let $X_0\subset \P^3$ be one of Shioda surfaces over a finite field $\F_p$ (\cite{Sh}). This means
that $X_0$ is smooth, unirational, but has large degree. Thus   in particular $H^2(X_0, \sO_{X_0})\neq 0$. Since $CH_0(X\times_{\F_p} \overline{\F_p(X)})=\Q$, \cite{E}, Lemma 2.1 implies that
$H^i(X_0, W\sO_{X_0, K})=0$  for all $i>0$.
Let $X\subset \P^4$ be the cone over $X_0$, and $\sigma: Y\to X$ be the blow up
of the vertex $v$. Then one has $R^2\sigma_*\sO_Y\surj H^2(X_0, \sO_{X_0})\neq 0$, thus the singularity is not rational in the sense of Remark \ref{rmk2.5}, a fortiori in the sense of Kempf. However the singularity is Witt-rational. Indeed, by \cite{BBE}, Theorem 2.4, applied to $\sI=$ maximal ideal of the vertex and $\sI'=\sO_Y(-X_0)$, one has $R^i\sigma_*\sI_K=0$ for $i\ge 1$ thus
\ga{2.8}{R^i\sigma_*W\sO_{Y,K}=\begin{cases} W\sO_{X,K} & i=0\\
H^i(X_0, W\sO_{X_0,K})=0 & i\ge 1 \end{cases}}
where for $i\ge 1$ it means the skyscraper sheaf with this value supported in $v$, showing that $X$ has Witt-rational singularities.
\section{Homological Chow group and alterations}
\subsection{Chow groups}
As explained in \cite{BEL}, sections 1 and 4, the right motivic condition
for $X$ projective over a perfect field $k$ of characteristic $p>0$ which forces $H^i(X/K)^{<1}=0$ for all $i\ge 1$
is as follows. Let $X\subset \P^N$ be a projective embedding, $\sigma: \P\to \P^N$ be an alteration so that $Y:=\sigma^{-1}(X)$ is a normal crossings
divisor. Then the graph of the alteration $\Gamma_\sigma\in H^{2N}(\P\times (\P^n\setminus X), Y\times (\P^N\setminus X), N)\otimes \Q$ should be supported along a divisor $A\subset \P$ which is in good position relatively to $Y$. Hypersurfaces of degree $\le N$ in $\P^N$ fulfill this condition (\cite{BEL},
Theorem 1.1) but in general, given $X$ this is difficult to compute as $X$ does not have a privileged embedding.

If we do not give ourselves an embedding but consider conditions defined
directly  on $X$, then a very strong one is to assume that motivic cohomology fulfills $H^{2n}(X, X_{{\rm sing}},n)\otimes \Q=\Q$. Here one has to pay attention. Relative motivic cohomology is defined when $X_{{\rm sing}}\subset X$ is a normal crossings divisor (see \cite{Le}  Chapter 4, 2.2 and p. 209), so is here not appropriate. But we could consider geometric conditions which would imply this motivic condition if we had resolution of singularities in charateristic $p>0$. One such geometric condition is to assume that $X \setminus X_{{\rm sing}}$ is rationally connected in the sense that any two closed points can be linked by a connected chain of rational projective curves in $X \setminus X_{{\rm sing}}$. Indeed, take $X_0\subset \P^N$ smooth of degree $N+1$, and define $X\subset \P^{N+1}$ to be the cone over $X_0$. Then for $N$ large enough, $X_0$ is not rationally connected, $(X\setminus X_{{\rm sing}})$ is not rationally connected, yet $X$  fulfills the conditon $H^i(X/K)^{<1}=0$ for all $i\ge 1$
(\cite{BEL}, Theorem 1.1), thus over a finite field $\F_q$,
one has $|X(\F_q)|\equiv 1 $ modulo $q$ (which is Ax' theorem).

The purpose of section \ref{sec3.2} is to show that the triviality of the homological Chow group $CH_0$, which is bad to count points  as $CH$ does not act on cohomology, nevertheless has a  nontrivial consequence on the part of cohomology surviving in alterations.
\subsection{$CH$ and the image of the cohomology in the cohomology of an alteration}\label{sec3.2}
We now show Proposition \ref{prop1.2} (see Introduction).
\begin{proof}[Proof of Proposition \ref{prop1.2}]
We apply Bloch's decomposition of cycles  in \cite{B}, Appendix to lecture 1. Since it is written for $X$ smooth proper, and applied to the diagonal, and we need it for $X$ projective not smooth, applied to the graph of the alteration, let us redo the argument.

First we observe that the kernel of $CH_0(X\times_k k(Y))\to CH_0(X\times_k K)$ is torsion, for $K=\overline{k(X)}$. Indeed, if a cycle $\xi\in CH_0(X\times_k k(Y))$ fulfills $\xi=\sum {\rm div}(f_W)$, with
$W$ curve on $X\times_k K$ and $f_W\in K(W)^\times$, then $(W,f_W)$
is defined over a finite extension $K_0$ of $k(Y)$, of degree $d$ say, so $d\xi=\sum_W {\rm div}{\rm Nm}_{(K_0/k(Y))}(f_W)$.

Then one has
\ga{3.1}{CH_0(X\times_k k(Y))=\varinjlim_{\emptyset\neq V {\rm open} \subset Y} CH_n(X\times V), \ n={\rm dim}(Y).}
Indeed, if $\xi\in CH_0(X\times_k k(Y))$, it is defined over some $\emptyset\neq V\subset Y$
and similarly if $(W, f_W)$ is a relation, it is defined over some open $V\neq \emptyset$ as well.

We consider the class of ${\rm Spec}(k(X)\times _{k(X)} k(Y))\cong {\rm Spec}(k(Y))$ in $CH_0(X\times_k k(Y))$. In view of \eqref{3.2}, this is the restriction to $X\times_k k(Y)$ of the class of the graph of the alteration $\Delta_\sigma\subset X\times Y$. The assumption
$CH_0(X\times_k K)=\Q$ implies that there is a 0-cycle $\xi\in CH_0(X)$ and $N\in \Z\setminus \{0\}$ so that $N\Delta_\sigma\equiv \xi \in CH_0(X\times_k k(Y))$. Thus we conclude by \eqref{3.2} that there is also a dimension $n$ cycle $\Gamma\subset X\times_k Y$, and a divisor $A\subset Y$ with
\ga{3.2}{N\Delta_\sigma\equiv \xi\times Y + \Gamma \in CH_n(X\times_k Y), \ \Gamma\subset X\times A.}
Let $X\subset \P^{n+c}=:\P$ be an embedding of $X$ (recall $n={\rm dim}(X)$).
By \cite{Fu}, Chapter 19, one has a homomorphism
\ga{3.3}{CH_n(X\times Y)\to H^{2(n+c)}_{X\times Y}(\P \times Y)(n+c).}
Since this is an important point in the argument, and the reference \cite{Fu} is written for Borel-Moore homology over $\C$, let us rewrite the argument.
Indeed, for an irreducible variety $W$ of dimension $n$ in $X\times Y$, one has by purity a cycle class in $H^{2(n+c)}_W(\P\times Y)(n+c)$, which maps to
$H^{2(n+c)}_{X\times Y}(\P \times Y)(n+c)$. And if $f:\sW\to \P^1$ is a rational map, yielding
\ga{3.4}{
\begin{CD}
\sW@>\subset >> X\times Y@>\subset >> \P\times Y\\
@VfVV\\
\P^1
\end{CD}
}
then if $\sW_t$ is a fiber of $f$ over $t$, of dimension $n$, one has that its cycle class
in $H^{2(n+c)}_{X\times Y}(\P \times Y)(n+c)$ is the image of the class of $t\in H^2(\sW)(1)$
via
\ga{3.5}{ H^2(\sW)(1) \xrightarrow{{\rm Gysin}}H^{2(n+c)}_{\sW}(\P \times Y)(n+c)\to
H^{2(n+c)}_{X\times Y}(\P \times Y)(n+c).}
So the difference of the classes of two fibers maps to 0.
Now one has an action
\ga{3.6}{H^i(X)\times H^{2(n+c) -i}_X(\P) \otimes H^i(Y)\to
 H^i(Y)(-n-c)}
via the duality on the first component. We apply this to the class
\ga{3.7}{\Delta_\sigma \in H^{2(n+c)}_{X\times Y}(\P\times Y)(n+c)= \\
 \ \sum_i [\Delta_i]\in
\oplus_i H^{2(n+c)-i}_X(\P)\otimes H^i(Y)(n+c),\notag}
where the decomposition on the second line is the K\"unneth decomposition.
One has
\ga{3.8}{[N\Delta_\sigma]=[\xi\times X]+ [\Gamma]}
Now
\ga{3.9}{[\xi\times X]_i\in {\rm Im} \ H^{2(n+c)-i}_{|\xi|}(\P)\otimes H^i(Y)(n+c)}
thus by purity
\ga{3.10}{[\xi\times X]_i=0 \ {\rm for}  \ i\neq 0.}
On the other hand
\ga{3.11}{[\Gamma]_i \in {\rm Im} H^{2(n+c)-i}_X(\P)\otimes H^i_A(Y)(n+c).}
Thus
\ga{3.12}{\sigma^*: H^i(X)\to N^1H^i(Y) \ \forall i\ge 1.}
This finishes the proof.
 \end{proof}
We now prove the main theorem of this note (see Theorem \ref{thm1.1} of the Introduction).
\begin{proof}[Proof of Theorem \ref{thm1.1}]
As $X$ has Witt-rational singularities, the map $\sigma^*: H^i(X/K)^{<1}\to H^i(Y/K)^{<1}$ is injective by Proposition \ref{prop2.4}. On the other hand, the Chow group condition implies by Proposition \ref{prop1.2}
that $\sigma^*H^i(X/K)\subset N^1H^i(Y/K)$.
But \cite{EFano}, Lemma 2.1 implies that $N^1H^i(Y/K)\cap H^i(Y/K)^{<1}=0$. Thus we conclude
\ga{3.13}{H^i(X/K)^{<1}=0.}
We then apply the Lefschetz trace formula for rigid cohomology
(\cite{ELS}, Th\'eor\`eme II)
\ga{3.14}{|X(\F_q)|=1+\sum_{i>0}(-1)^i {\rm Tr}\ F|H^i(X/K)}
where $F$ is the geometric Frobenius of $\F_q$. This finishes the proof.
\end{proof}

\subsection{Koll\'ar's example} \label{sec3.3}
J. Koll\'ar constructed an example of a rationally connected surface $X/\F_q$
without rational point. Since rational connectedness implies $CH_0(X\times \overline{\F_q(X)})=\Q$ ,  Theorem
\ref{thm1.1} implies  that the singularities are worse than Witt-rational. Here is the example.
Let $Y'=C\times \P^1$, where $C$ is a genus $\ge 2$ curve with $C(\F_q)=\emptyset$, let $ \ y=(\xi\times \alpha)\in Y', \ \sigma: Y\to Y'$ be the blow up of the closed point $y$ (which is not rational), and set $\tau^{-1}(C)=C'+E$ where $E$ is the exceptional locus and $C'$ is the proper transform of $C\times \alpha$. Then one has
\ga{3.16}{(C')^2=C^2-{\rm deg}(\xi\times \alpha)=-{\rm deg}(\xi\times \alpha),}
which one computes as when the ground field is algebraically closed by writing
\ga{3.17}{(C'+E)\cdot C'= C^2=0.}
Thus we conclude that
\ga{3.18}{(C')^2<0.}
Thus if $k=\F_q$,  by Artin theorem (\cite{Ar}, Theorem 2.9), there is  a contraction $\sigma: Y\to X$,
where  $X$ is a projective surface with $\sigma|_{Y\setminus C'}={\rm isomorphism}$ and $\sigma(C')=x\in X$. This yields the correspondence
\ga{3.19}{\begin{CD}
X\\
@A\sigma AA \\
Y@>\tau >> Y'
\end{CD}
}
If deg $\alpha >1$, then $X$ has no rational point and yet $X\times_k K$
is rationally connected for any algebraically closed field $K$ containing $k$.

We explain now  this example cohomologically.
We denote by $H(Z)$ rigid cohomology $H(Z/K)$. As $\tau$ is a blow up on a smooth surface, one has
\ga{3.20}{\tau^*: H^1(Y')\xrightarrow{\cong} H^1(Y)}
therefore
\ga{3.21}{H^1(Y)\xrightarrow{{\rm rest}} H^1(C')}
is injective and not surjective as $\alpha$ is not a rational point. For example, if we request $\alpha$ to be of degree 2, then
$H^1(C')/H^1(Y)\cong H^1(C)$. Consequently
\ga{3.22}{(H^1(C')/H^1(Y))^{<1}\neq 0.}
This yields
\ga{3.23}{ H^1_c(Y\setminus C')=H^1(X)=0\\
0\to   H^1(C')/H^1(Y)\to H^2_c(Y\setminus C')=H^2(X) \xrightarrow{\sigma^*} H^2(Y)\notag}
and violates Proposition \ref{2.4}. Therefore the singularities of $X$ are not Witt-rational.

\bibliographystyle{plain}

\begin{thebibliography}{99}
\bibitem{Ar} Artin, M.: Some numerical criteria for contractibility of curves on algebraic surfaces, Am. J. of Mathematics, {\bf 84} (1962), 485-496.
\bibitem{BBE} Berthelot, P., Bloch, S., Esnault, H.: On Witt vector cohomology
for singular varieties, preprint 2005, 42 pages.
\bibitem{B} Bloch, S.: Lecture on Algebraic cycles, Duke
University Mathematics Series, IV, (1980).
\bibitem{BEL} Bloch, S., Esnault, H.,  Levine, M.: Decomposition of the
diagonal and eigenvalues of Frobenius for Fano hypersurfaces,
 Am. J. of  Mathematics,  {\bf 127} no1 (2005), 193-207.
\bibitem{deJ} de Jong, A. J.: Smoothness, semi-stability and alterations, Publ. Math. IHES {\bf 83} (1996), 51-93.
\bibitem{DeInt} Deligne, P.: Th\'eor\`eme d'int\'egralit\'e, Appendix  to
 Katz, N.: Le niveau de la cohomologie des intersections compl\`etes,
Expos\'e XXI
in SGA 7, Lect. Notes Math. vol. {\bf 340}, 363-400,
Berlin Heidelberg New York Springer 1973.
\bibitem{E} Esnault, H.: Hodge type of subvarieties of $\P^n$ of small
degrees. Math. Ann.  {\bf 288}  (1990),  no. 3, 549-551.
\bibitem{EFano} Esnault, H.: Varieties over a finite field with trivial Chow group of 0-cycles
have a rational point, Invent. math. {\bf 151} (2003), 187-191.
\bibitem{Euneq} Esnault, H.:  Deligne's integrality theorem in unequal characteristic and rational
points over finite fields, preprint 2004, 10 pages, appears in the Annals
of Mathematics.
\bibitem{ELS} \'Etesse, J.-Y.; Le Stum, B.: Fonctions $L$ associ\'ees
aux $F$-isocristaux surconvergents. I, Math. Ann. {\bf 296} (1993),
557--576.
\bibitem{Fu} Fulton, W.: Intersection Theory, Ergebnisse der Mathematik und ihrer Grenzbegiete, {\bf 2} Springer Verlag 1984.
\bibitem{Gr} Grothendieck, A.: Formule de Lefschetz et rationalit\'e
des fonctions $L$. S\'eminaire Bourbaki {\bf 279}, 17-i\`eme
ann\'ee (1964/1965), 1-15.
\bibitem{Ke} Kempf, G.:
Varieties with rational  singularities,
The Lefschetz centennial conference, Part I (Mexico City, 1984),
  Contemp. Math., {\bf 58}, 179-182,
Amer. Math. Soc., Providence, RI, 1986.
\bibitem{Ko} Kov\'acs, S.: A characterization of rational singularities, Duke Math. J. {\bf 102} 2 (2000), 187-191.
\bibitem{Le} Levine, M.: Mixed Motives, Mathematical Surveys and Monographs
{\bf 57} (1998), American Mathematical Society.
\bibitem{Sh} Shioda, T.:
An example of unirational surfaces in characteristic $p$,  Math. Ann.  {\bf 211}  (1974), 233-236.
\end{thebibliography}

\renewcommand\refname{References}

\end{document}